\documentclass[a4,11pt]{amsart}
\usepackage{amssymb,amsmath,amsthm,amsfonts}
\theoremstyle{plain}
\newtheorem{thm}[subsection]{Theorem}

\newtheorem{lemma}[subsection]{Lemma}
\newtheorem{cor}[subsection]{Corollary}

\theoremstyle{definition}

\newtheorem{rmk}[subsection]{Remark}
\newtheorem{eg}[subsection]{Example}

\def\<{\langle}
\def\>{\rangle}

\newtheorem*{oldproof}{Proof}
\renewenvironment{proof}[1][{}]{\begin{oldproof}[#1]}{\qed\end{oldproof}}

\begin{document}
\title[ finite non-commutative algebras]{maximum number of solutions of
$x^q =x $ in a finite non-commutative algebra}
\author{Vineeth Reddy Chintala}
\address{Department of Mathematics, Tata Institute of Fundamental Research, Mumbai,
India.} \email{vineethreddy90@gmail.com}
\date{\today}
\subjclass[2000]{} 

\begin{abstract}

 We obtain a bound on the number of solutions of $x^q=x$
in a finite noncommutative algebra over a field with $q$ elements. 
Furthermore, we completely characterize those rings for which this maximum number is attained.

\end{abstract}
\maketitle

\section*{\textbf{Introduction}}

\indent If every element of a ring is an idempotent, then the ring is commutative. What is the highest proportion of 
idempotents possible in a finite noncommutative ring?
More generally, Jacobson proved that if all the elements of a ring satisfy
the equation $x^{n(x)} = x$ for some $n(x)>1$ depending on $x$, then the ring is commutative (\cite{TYL} Theorem 12.10). 
And similarly we can ask about the number of solutions of $x^n = x$ in a finite noncommutative ring.

\indent In this article we work with noncommutative algebras over a finite field with $q$ elements.
We consider the set $I= \{x | x^q - x \in C \}$, where $C$ is the center and get a bound on the percentage of elements 
that fall into $I$. 
The pivotal result is Lemma \ref{main}, which states a condition for two elements to commute with each other.
It is used in measuring the maximum density of $I$ in the rings under consideration. \\
\indent In particular, we get a bound on the maximum number of solutions of $x^q =x$ in a finite algebra over $F_q$. 
We also characterize the algebras in which this maximum is attained. 

\indent  The last part contains a few observations about idempotents in a ring. 
We prove that the number of idempotents in any finite noncommutative ring cannot exceed
three-fourths of its cardinality. The article is for the most part elementary. \\

\noindent
An \emph{idempotent} $x$ in a ring $R$ is an element such that $x^2=x$ and the ring $R$ is said to be Boolean
if each element is an idempotent. We will denote the center of the ring by $C$. For each element $b$ of the ring, let $C_b$
denote the set of all elements of the ring which commute with $b$. A ring is said to be indecomposable if it is not 
a direct product of two rings. All rings considered in this article are associative with unity. 

\vskip 5mm

\section {\textbf{Commutativity of two elements : $ab=ba$?}}

\vskip 3mm

\indent Let $F_q$ be a finite field with $q$ elements $f_1, f_2....,f_q$.
Observe that if $$s := \sum_{i=1}^{q} (f_i) ^r$$ then $s$ is either $-1$ or
$0$ according as to whether $q-1$ divides $r$ or not.
Indeed, this follows simply by multiplying $s$ by $f^r$ where $f$ is
a generator of the cyclic group $F_q^{\ast}$. Since $f^r s =
s$ in $F_q$, we have $s=0$ unless $f^r=1$ which happens only when $(q-1)|r$.

\vskip 5mm

\begin{lemma}\label{main}
  Let $R$ be an algebra over $F_q$ and $P$ be a polynomial of degree $q$ with coefficients in $F_q$.
Let $a , b$ be any two elements of $R$ such that $b^q - b$ lies in the center of $R$. If  $ \sum_{i=1}^{q} P(a + f_ib) \in C_b$,
then $ab = ba$.
\end{lemma}
\vskip 3mm

\begin{proof}
Since we are working over a field, we can take the polynomial to be monic.
Let $$P(x) = x^q+r_{q-1}x^{q-1}+.....+r_0$$ where $r_i$ are elements of $F_q$.
Note that
$$P(a + tb) = P(a) +tc_1+t^2c_2+.....+t^{q-1}c_{q-1}+t^qb^q$$ where $c_i$'s are certain elements of $R$ that are independent of $t$.
Summing over $t$ from $f_1$ to $f_q$ gives (by the observation about $\sum_{i=1}^{q} (f_i) ^r$)
that $c_{q-1} \in C_b$. Note that $$c_{q-1} = r_{q-1}b^{q-1} + b^{q-1}a + b^{q-2}ab+ \cdots
+ab^{q-1}.$$ This equation gives on pre-multiplying and post-multiplying by $b$,
the respective equations $$bc_{q-1} = r_{q-1}b^{q} + b^q a + b^{q-1}ab+ \cdots +bab^{q-1},$$
and $$c_{q-1}b = r_{q-1}b^{q} +b^{q-1}ab+ \cdots +bab^{q-1}+ab^q.$$ Subtracting, we get
$ab^q=b^qa$. Since $b^q - b$ lies in the center, this implies that $a$ and $b$ commute with each other.
\end{proof}

\vskip 5mm

\begin{cor}\label{cor1}
Let $R$ be an algebra over $F_q$ and $b \in R$ such that $b^q - b \in C$. 
If $ab \neq ba$, then $(a +fb)^q - (a+fb) \notin C_b$ for some $f \in F_q$.
\end{cor}
\vskip 5mm

\begin{cor}
Let $R$ be an algebra over $F_q$ and $b \in R$ such that $b^q - b \in C$. 
If $x^q -x  \in C_b$ for all elements $x \in R$, then $b$ lies in the center of the ring.
\end{cor}
\vskip 5mm

\begin{rmk}
Let $I = \{x| x^q - x \in C \}$, where $C$ is the center of $R$.
If $I$ is an additive subspace of $R$, then the elements of $I$ commute with each other.
In particular, if $I = R$ then $R$ is a commutative algebra. This gives us an elementary proof
of a special case of Jacobson's theorem : if $x^q =x$ is satisfied throughout in an algebra
over $F_q$, then the algebra is commutative.
\end{rmk}
\vskip 5mm

\begin{thm}\label{thm1}
Let $R$ be a finite noncommutative algebra of $q^n$ elements over $F_q$.
Let $I = \{x| x^q - x \in C \}$ and $r =\frac{|I|}{|R|}$. 
Then $r \leq \frac{q^2 - q + 1}{q^2}$. Further, if  $r =  \frac{q^2 - q
+ 1}{q^2}$ then $R$ is generated as an algebra by the center $C$ and
two non-commuting elements $b,c \in I$ such that $C\subseteq I$, $|C| =q^{n-2}$
and $|C_b|= |C_c| = q^{n-1}$.

\end{thm}
\vskip 3mm

\begin{proof}

 If $I\subseteq C$, then $|I|\leq|C|\leq q^{n-2}$.
Let $b\in I, b \not\in C$. Consider the additive subgroup  $B = \{ub : u \in  F_q\} $. 
Then  $R =\bigcup_{i=1}^{q^{n-1}}(a_i + B)$, where $(a_i+B)$ is a coset of $B$.

Since $C_b$ is a proper subalgebra of $R$, we have $|C_b| \leq q^{n-1}$. So there are at least $q^n - q^{n-1}$ 
elements which do not commute with $b$. These elements fall into $q^{n-1} - q^{n-2}$ cosets of $B$.
By Corollary \ref{cor1}, each of these cosets contains an element which does not belong to $I$. 
Hence there are at least $q^{n-1} - q^{n-2}$
elements lying outside $I$. Therefore $r \leq \frac{q^2 - q + 1}{q^2}$.\\

Suppose $r = \frac{q^2 - q + 1}{q^2}$. It follows from the proof above that $C_b\subseteq I$ and $|C_b|= q^{n-1}$.
Let $c\in I, c \not\in C_b$. Similarly, we can infer that $C_c\subseteq I$ and
$|C_c|= q^{n-1}$. Clearly, we have $R = B + C_c$.
Note that $C_b$ and $C_c$ are commutative subrings of $R$ and $C_b\cap C_c = C$.

\indent We will now prove that $|C|= q^{n-2}$. Let $ t= |C_c|/|C|$; then $C_c =\bigcup_{i=1}^{t}(c_i + C)$,
where $(c_i+ C)$ is an additive coset of $C$ relative to $C_c$.
Suppose $c_i\notin C$. It follows from Corollary \ref{cor1} that $c_i +B$ contains an element which does not belong to $I$.
Therefore the set $c_i + B + C$ has at least $|C|$ elements lying outside $I$.
Note that $R = \bigcup_{i=1}^{t}(c_i + B + C)$; counting the elements of $R\backslash I$ will give us
$$(t-1)|C|\leq |R\backslash I| = (q-1)q^{n-2}$$
Since $ t|C| = |C_c|= q^{n-1}$, it follows that $|C| = q^{n-2}$.
\end{proof}
 
\vskip 5mm

\section{\textbf{Solutions of $x^q=x$ in finite noncommutative algebras}}    
\vskip 2mm

\emph{Let $I_q$ be the set of solutions of $x^q=x$ in $R$ and $r_q = \frac{|I_q|}{|R|}$.} 
The following result follows directly from Theorem \ref{thm1}.

\vskip 3mm

\begin{thm}\label{thm2}
Let $R$ be a finite noncommutative algebra of $q^n$ elements over $F_q$.
Then $r_q \leq \frac{q^2 - q + 1}{q^2}$. Further, if  $r_q =  \frac{q^2 - q
+ 1}{q^2}$ then $R$ is generated as an algebra by the center $C$ and
two non-commuting elements $b,c \in I_q$ such that $C\subseteq I_q$, $|C| =q^{n-2}$
and $|C_b|= |C_c| = q^{n-1}$.

\end{thm}

\vskip 5mm

\begin{cor}\label{cor3}
Let $p \geq 2$ be a prime and $R$ a finite noncommutative ring of $p^n$ elements.
Then $r_p \leq \frac{p^2 - p + 1}{p^2}$. Further, if  $r_p =  \frac{p^2 - p
+ 1}{p^2}$ then $R$ is generated as an algebra by the center $C$ and
two non-commuting elements $b,c \in I_p$ such that $C\subseteq I_p$, $|C| =p^{n-2}$
and $|C_b|= |C_c| = p^{n-1}$.
\end{cor}
\vskip 2mm

\begin{proof}
 Suppose  $\textit{char}$ $R = p^k$. If $k = 1$, then $R$ is a noncommutative algebra over $F_p$ and we are down 
to a special case of Theorem \ref{thm2}. \\
\indent Suppose $k>1$. Then we have
 $(x + p^ {k-1}) ^ p = x^p.$
Therefore $I_p \cap (p^{k-1} + I_p)  = \emptyset$, implying that $I_p \leq \frac{|R|}{2}$ in this case.  
\end{proof}

\vskip 3mm

\begin{cor}\label{cor2}
Let $R$ be a finite algebra over $F_q$. 
If $r_q > \frac{q^2 - q + 1}{q^2}$, then $R\simeq F_q^n$. 
Such rings will be called \emph{$q$-rings.} (Here $q$ denotes the cardinality of the underlying field.)
\end{cor}
\vskip 2mm

\begin{proof}
It follows from Theorem \ref{thm2} that $R$ is commutative.
In that case $I_q$ is an additive subgroup of $R$. Therefore $|I_q|$ divides $|R|$ and $r_q = 1$. Now $R$ is finite, 
commutative and has no nilpotent elements. Therefore $R$ is semisimple and is isomorphic to a direct product of fields.
These fields are algebras over $F_q$ and have at most $q$ elements since the elements of the field satisfy $x^q =x$.
Therefore $R\simeq F_q^{n}$. \emph{Note that $R$ has $2^n$ idempotents.}
\end{proof}
\vskip 5mm

Here is an example where $r_q =  \frac{q^2 - q + 1}{q^2}$.

\begin{eg}\label{egp}
Consider the algebra $S$ of cardinality $q^3$ generated by $F_q$ and
symbols $e_1,e_2$ with the following operations :
$e_ie_j = e_i$, $1\leq i,j \leq 2$.

Note that $(e_1 - e_2)^2 = 0$. It follows that for every $f_1, f_2 \in F_q$ and $m\geq1$,
$$[f_1(e_1 - e_2) + f_2]^m = f_2^m + m{f_2}^{m-1}f_1(e_1 - e_2).$$
Also,
\begin{align*}
  (f_1e_1+ f_2e_2)^m &= f_1e_1(f_1e_1+ f_2e_2)^{m-1} +  f_2e_2(f_1e_1+ f_2e_2)^{m-1}\\
                  &= f_1e_1(f_1+f_2)^{m-1} + f_2e_2(f_1+f_2)^{m-1}.
\end{align*}
Therefore $(f_1e_1+ f_2e_2)^q \neq (f_1e_1+ f_2e_2)$ if and only if $f_1 + f_2 = 0$.
And so $S$ has exactly $q(q^2 - q +1)$ solutions of the equation $x^q = x$.
By taking the product of this ring with a $q$-ring of cardinality $q^n$, we get a ring of cardinality $q^{n+3}$ satisfying
the equality in Theorem \ref{thm2}.
\end{eg}

\vskip 5mm

\begin{rmk}

Let $i_p$ denote the number of solutions of $x^p = x$ in a ring. A direct consequence of the above discussion is 
that there is a sequence of finite noncommutative
rings $\{R_p$: $p$ is prime$\}$ such that $\frac{i_p}{|R_p|}$ tends to $1$.
\end{rmk}

\vskip 5mm

\begin{thm}
Let $R$ be a finite indecomposable algebra over $F_q$. If $r_q = \frac{q^2 - q + 1}{q^2}$,
then $R \simeq S$, where $S$ is the ring defined in Example \ref{egp}.
\end{thm}
\vskip 2mm

\begin{proof}
 If $R$ is commutative, then $I_q$ is a subring of $R$ and
 $|I_q|$ divides $|R|$. Thus the hypothesis
implies that $R$ is noncommutative. 

   Suppose $R$ has $q^n$ elements. If $n < 3$, then $R$ is commutative. Therefore $n\geq3$. It follows from
Theorem \ref{thm2}  that the center of $R$ is a $q$-ring of $q^{n-2}$ elements
(with $2^{n-2}$ idempotents). If $n>3$, then $R$ has a nontrivial central idempotent
and so it can be written as a product of two rings, contradicting the indecomposability of $R$. Therefore $n=3$.
It is enough if there is a unique $R$ satisfying the conditions. 
\end{proof}
\vskip 3mm

The following result is presumably well known.
\vskip 3mm

\begin{thm}
Let $R$ be a noncommutative algebra with unity over a field $F$ (possibly infinite).
If $\textit{dim}$ $R= 3$ as a vector space over $F$, then $R$ is unique up to isomorphism.
\end{thm}
\vskip 2mm

\begin{proof}
By a straightforward case-by-case analysis, it can be shown that the Jacobson radical $J$ is one-dimensional
and $R/J\simeq F\times F$; then the multiplication table of $R$ is easily determined.
\end{proof}

\vskip 5mm

\section{\textbf{Idempotents in finite rings}}

\vskip 3mm
\emph{Let $I_2$ denote the set of idempotents in a finite ring $R$ and $i = |I_2|$.} 

\vskip 3mm

\begin{thm}
 If $R$ is noncommutative, then $ i \leq \frac{3|R|}{4}$. Further, if $i = \frac{3|R|}{4}$ then
$R$ is generated as a ring by a Boolean ring $C$ contained in the center of $R$ and two noncommuting idempotents $b,c$
such that $|C| = \frac{|R|}{4}$.
\end{thm}

\vskip 2mm

\begin{proof}
If $R = R_1\times R_2$ then $ \frac{i}{|R|} = \frac{i_1}{|R_1|}\frac{i_2}{|R_2|}$, where $i_j$ is the 
number of idempotents in $R_j$. Since $\frac{i_j}{|R_j|} \leq 1$, it is enough to consider the case when
$R$ is indecomposable. 

\indent Then  $|R| = p^n$, where $p$ is a prime. If $p >2$ then at most one of $x$, $-x$ is an idempotent unless $x=0$. 
Hence $i \leq \frac{|R|}{2} + 1$ when $p>2$. When $p =2$, the statement is reduced to a special case of 
Corollary \ref{cor3}. 
\end{proof}

\vskip 3mm

\begin{cor}\label{cor4}
 Let R be a finite ring with $i >\frac{3|R|}{4}$. Then R is Boolean.
\end{cor}

\vskip 2mm

\begin{proof}
 Since $i >\frac{3|R|}{4}$, $R$ is commutative. Now $R = R_1\times...\times R_n$
where each $R_t$ is an indecomposable ring. Since each $R_t$ has exactly two idempotents, $R$ has $2^n$ idempotents. 
So $2^n > \frac{3}{4}|R| = \frac{3}{4}|R_1|...|R_n|\geq \frac{3}{4}.2^{n-1}|R_i| $. This implies that $2> \frac{3}{4}|R_t|$.
Therefore $|R_t| = 2$ and so $R = F_2{^n}$ is boolean.
\end{proof}

\vskip 5mm

Arguing similarly we can prove the following two theorems.

\vskip 3mm 

\begin{thm}
 Let $p \geq 2 $ be a prime and $R$ a commutative ring such that $p$ divides $|R|$. Then $\frac{i}{|R|} \leq \frac{2}{p}$,
equality occurring if and only if $R\simeq F_p\times F_2^{n-1}$.
\end{thm}
\vskip 2mm

\begin{proof}
  We can write $R = R_1\times...\times R_n$ such that $p$ divides $|R_1|$ and each $R_t$ is an indecomposable ring.
So $2^n =\frac{i}{|R|}|R_1|...|R_n|\geq \frac{i}{|R|}2^{n-1}|R_1|$.
Therefore $\frac{i}{|R|}\leq \frac{2}{|R_1|}\leq\frac{2}{p}$, equality occurring if and only if $|R_1| = p$
and $|R_t| = 2 (t\neq1)$, i.e. $R\simeq F_p\times F_2^{n-1}$.
\end{proof}
\vskip 5mm

\begin{thm}
For a finite ring $R$, the following are equivalent.
\\
1. $R$ is boolean.
\\
2. $i > \frac{3}{4}|R|$.
\\
3. $i_c > \frac{2}{3}|R|$, where $i_c$ is the number of central idempotents in $R$.
\\
4. $i_c > \frac{1}{2}|R|$ and $3\nmid |R|$.
\\
Moreover, we have $\frac{|R|}{2} < i_c \leq \frac{2|R|}{3}$ if and only if $R\simeq F_3\times F_2^{n-1}$.
\end{thm}

\vskip 2mm

\begin{proof}
 Each of the four conditions imply that $R$ is commutative. Indeed this follows from Corollary \ref{cor4} and
the fact that $|C|$ divides $|R|$, where $C$ is the center of the ring. The proof is similar to that of 
Corollary \ref{cor4}.
\end{proof}

\vskip 5mm

\noindent {\bf Acknowledgements.} \vskip 5mm
I thank Prof. B. Sury for his suggestions which improved the readability of the article. The referee's comments were 
also helpful.

\end{document}